\documentclass[11pt]{elsarticle}
\usepackage{lineno,hyperref}
\usepackage{optidef}
\usepackage{graphicx}
\usepackage{float}
\usepackage{graphicx}
\usepackage{subfig}
\usepackage{mathrsfs}
\usepackage{float}
\usepackage{amsmath}
\usepackage{amsfonts}
\usepackage{booktabs}
\usepackage{natbib}
\usepackage{bbm}
\usepackage{algorithm}
\usepackage[noend]{algpseudocode}
\usepackage{multirow}
\algdef{SE}[DOWHILE]{Do}{doWhile}{\algorithmicdo}[1]{\algorithmicwhile\ #1}
\usepackage[letterpaper, portrait, margin=1in]{geometry}
\usepackage{booktabs}
\usepackage{tabularx}
\usepackage{xcolor}
\hypersetup{
    colorlinks,
    linkcolor={red!50!red},
    citecolor={blue!50!blue},
    urlcolor={blue!80!blue}
}
\usepackage{optidef}

\usepackage{amsthm}

\theoremstyle{definition}

\usepackage{newtxtext} 

\makeatletter
\newcommand*{\centerfloat}{%
  \parindent \z@
  \leftskip \z@ \@plus 1fil \@minus \marginparwidth
  \rightskip \leftskip
  \parfillskip \z@skip}
\makeatother

\usepackage{amssymb,amsmath,caption}
\usepackage[hang,flushmargin]{footmisc}
\newcommand{\algorithmfootnote}[2][\footnotesize]{
  \let\old@algocf@finish\@algocf@finish
  \def\@algocf@finish{\old@algocf@finish
    \leavevmode\rlap{\begin{minipage}{\linewidth}
    #1#2
    \end{minipage}}%
  }%
}

\usepackage{setspace}
\usepackage{natbib}
\usepackage{algorithm}
\usepackage[noend]{algpseudocode}
\algdef{SE}[DOWHILE]{Do}{doWhile}{\algorithmicdo}[1]{\algorithmicwhile\ #1}

\usepackage{appendix}

\usepackage{array}
\usepackage{threeparttable}
\usepackage{multirow}
\usepackage{rotating}
\usepackage{makecell}

\journal{.}

\biboptions{authoryear}

\begin{document}
\begin{frontmatter}

\title{Self-Calibrated Transit Service Monitoring Using Automated Collected Data}


\author{Hongyu Guo}
\address[label1]{Hangzhou Basis International School, Hangzhou, China, 310020}
\begin{abstract}
This paper proposes a self-calibrated transit service monitoring framework that aims to obtain the performance of a transit system using automated collected data. We first introduce an event-based transit simulation model, which allows the detailed simulation of passenger travel behavior in a transit system, including boarding, alighting, and transfer walking. To estimate passenger path choices, we assume the path choices can be modeled using a C-logit model, and propose a simulation-based optimization model to estimate the path choice parameters based on automated fare collection and automated vehicle location data. The path choices can be estimated on a daily basis, which enables the simulation model to adapt to dynamic passenger behavior changes, and output more accurate network performance indicators for regular service monitoring such as train load, passenger travel time, and crowding at platforms. The proposed system eliminates the need for conventional monitoring equipment such as cameras at platforms and scaling/weighing systems on trains. The Hong Kong Mass Transit Railway (MTR) system is used as the case study. Results show that the model can well estimate the path choice behavior of passengers in the system. The output passenger exit flows are closer to the actual one compared to the two benchmark models (shortest path and uniform path choice). 
\end{abstract}

\end{frontmatter}


\section{Introduction}\label{intro}
The global trend towards eco-friendly and effective transportation solutions has sparked a significant increase in the growth of public transit networks, particularly in China, Europe, and Japan. According to a study by \citet{metro_2021}, China boasts two of the world's most extensive metro systems, situated in Shanghai and Beijing, covering a total distance of 795.4 kilometers and 709.9 kilometers, respectively. Remarkably, Shanghai Metro holds the title of the world's longest public transport system. The expansion of transit networks has yielded numerous favorable outcomes, such as improved connectivity, decreased traffic congestion, and a reduction in carbon emissions.

As public transit networks continue to expand, the demand for these systems has experienced a significant upswing. For instance, the Beijing metro system's daily passenger traffic has surged from 1.2 million in 2000 to a staggering 9.3 million in 2023 \citep{wang2014analysis, beijing_metro_2023}. This heightened demand presents several challenges to the smooth operation of public transit systems. Among these challenges, overcrowding stands out as a prominent issue, particularly during peak hours when passengers inundate the system within a short time frame. The congestion on trains and at stations can result in passenger discomfort, reduced operational efficiency, and potential safety concerns.

Evaluating network performance and service quality, such as obtaining information about platform wait times and train loads, is crucial for transportation agencies and operators. It helps them gain insights into system functionality, provide valuable information to passengers, and enhance operational strategies. Nevertheless, the conventional methods of monitoring network performance often entail substantial investments in hardware and infrastructure, such as the deployment of cameras and passenger counting systems on platforms \citep{saponara2016exploiting, velastin2020detecting, seidel2021napc} and the implementation of scaling and weighing systems on trains \citep{nielsen2014estimating}. These investments can be cost-prohibitive, especially for large-scale transit systems.

An alternative approach to assess network performance involves utilizing automatically collected data, notably automated fare collection (AFC) data and automated vehicle location (AVL) data \citep{mo2022toward}. In contrast to hardware-dependent network performance monitoring, these methods, relying on automatically gathered data, prove to be more cost-effective as they eliminate the need for additional equipment. AFC data encompasses passenger usage transactions derived from smart cards, including location and time information. Depending on the fare system, AFC data may encompass transactions for both rail and bus travel (e.g., the urban rail system in Chicago) or exclusively for rail travel (e.g., the metro system in Hong Kong). AFC systems can be categorized as open or closed. Open systems require passengers to tap in only when they enter the system (as exemplified by the urban rail system in Boston), while closed systems necessitate both tapping in and tapping out (e.g., the transit system in Seoul, Korea). Many systems adopt a hybrid approach, combining an open architecture on the bus side and a closed one on the subway side, as seen in London. In the case of a closed system, AFC data can offer precise origins and destinations of passengers' journeys, along with tap-in and tap-out times. AVL data, on the other hand, provides information about the real-time location of vehicles, with train locations obtained from the rail tracking system and bus locations acquired from the GPS units installed in the vehicles. Given the wealth of information contained within AFC and AVL data, they present abundant opportunities for analysis in various domains, including travel behavior, operational planning, and performance monitoring for public transit systems.

\begin{figure}[htb] 
\centering
\includegraphics[width=0.8\linewidth]{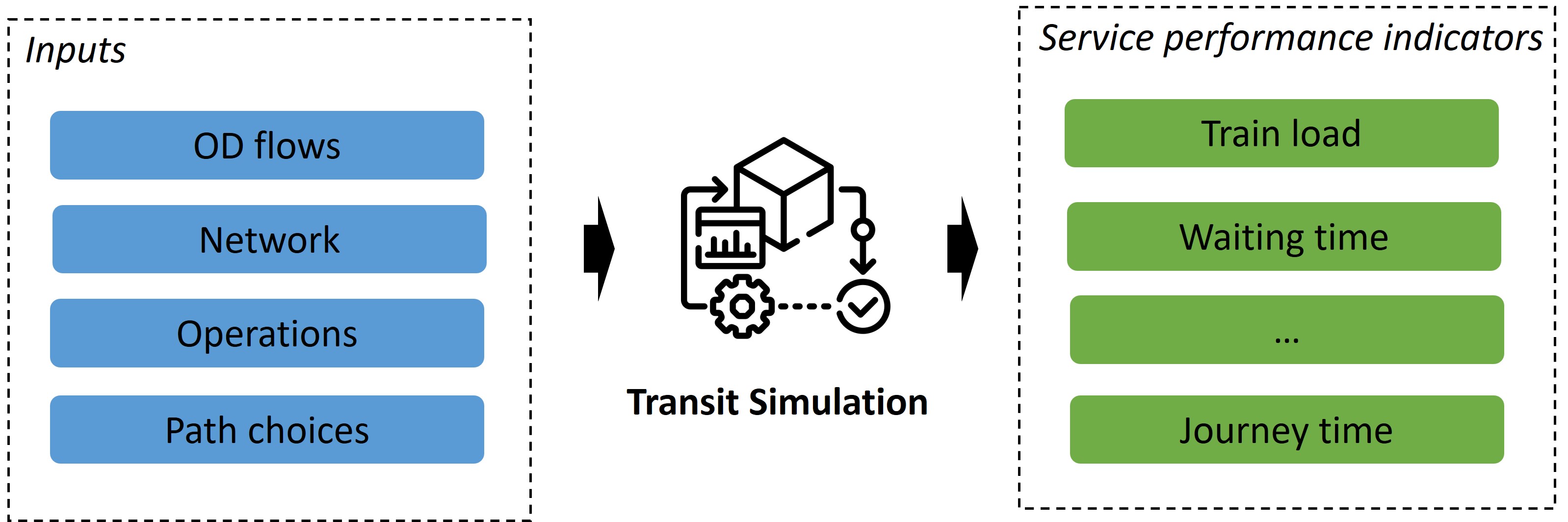}
\caption{Using transit simulation model for service performance monitoring}
\label{fig_intro_sim}
\end{figure}

In terms of performance monitoring, some key service level indicators can be directly extracted from AFC and AVL data, including metrics like vehicle kilometers, vehicle hours, travel time reliability, and Origin-Destination (OD) demand \citep{trepanier2009calculation, ma2014development}. However, accurately determining vehicle load and passenger waiting times poses a non-trivial challenge. Recently, various methods have been put forth to estimate passenger waiting times and the number of passengers left behind at stations, drawing upon AFC and AVL data \citep{ma2019estimation, zhu2017inferring}. Nevertheless, these approaches typically focus on estimating only one or two specific network performance indicators. To obtain a comprehensive set of performance indicators, it is often necessary to employ a transit simulation model \citep{mo2020capacity}. Such models offer detailed insights into service performance at various levels, including stations, lines, trains, and individual passengers. These simulations take into account factors such as OD flows, path choices, and operational specifics as input, providing a holistic view of the network's performance, as illustrated in Figure \ref{fig_intro_sim}.

There are two major tasks for simulation-based service performance monitoring. The first is to develop a transit simulation model that can mimic the passenger and vehicle interactions in a network. The second is to get all input data (including model parameters) for the simulation. OD flows, network, and operation details can be directly obtained from AFC and AVL data, but path choices are not directly observable. Therefore, the remaining challenge of the second step is to calibrate path choices in the transit simulation model,

This paper proposes a self-calibrated transit simulation model following the same framework as \citet{mo2020capacity}. Specifically, an event-based transit simulation model is proposed, which allows the detailed simulation of passenger travel behavior in a transit system, including boarding, alighting, and transfer walking. To estimate the path choices, we assume the path choices can be modeled using a C-logit model, and propose a simulation-based optimization model to estimate the path choice parameters based on AFC and AVL data. The path choices can be estimated on a daily basis, which enables the simulation model to adapt to dynamic passenger behavior changes, and output more accurate network performance indicators for regular service monitoring. The Hong Kong Mass Transit Railway (MTR) system is used as the case study. We generate synthetic smart card data with a set of hypothetical ``true`` path choice parameters. The reason for using synthetic data is to allow the comparison with ``ground truth''. Results show that the model can estimate the path choice behavior of passengers in the system. The output OD exit flows are closer to the actual one compared to the two benchmark models (shortest path and uniform path choice). 

The remainder of the paper is organized as follows: Section \ref{sec_liter} reviews recent literature on hardware-based and model-based transit service monitoring. Section \ref{sec_method} introduces the methodology, including the path choice model, transit simulation, and the path choice calibration model. Case studies are presented in Section \ref{ref_case} to validate the model's performance and demonstrate its functionality. Section \ref{sec_con_dis} concludes the paper and discusses future research directions.

\section{Literature review}\label{sec_liter}

\subsection{Hardware-based transit service monitoring}
The traditional way of monitoring transit services is through hardware settled up around the stations. For example, \citet{saponara2016exploiting} aimed to use the onboard closed-circuit television (CCTV) security system in trains to collect onboard passenger information. This information can be used to improve train services and safety, such as automatic train temperature control and smoking detection. \citet{nielsen2014estimating} introduced a new counting technique that utilizes the weighing systems found in modern trains.  The study compares the passenger counts obtained through this technique with manual counts and counts from an infrared system in trains in urban Copenhagen, showing results from the weighing system are better than those of the infrared equipment. \citet{seidel2021napc} presented a real-time Neural Automated Passenger Counting system (NAPC) that utilizes an end-to-end LSTM recurrent neural network. NAPC is designed to accurately count the number of passengers getting on and off a vehicle using 3D LiDAR videos captured by a top-down sensor during door opening phases. The results indicate the potential to create a privacy-friendly Automated Passenger Counting (APC) system with competitive counting performance, employing a deep-learning approach and minimizing unnecessary data collection in public spaces. \citet{velastin2020detecting} introduced a machine learning-based approach to locate, track, and count people in a public transport setting using data from standard video cameras. The study evaluates the performance of three state-of-the-art detectors and trackers, ultimately choosing D-SORT and Yolov3 models for people counting due to their computational efficiency.

\subsection{Model-based transit service monitoring using AFC and AVL data}
With the emergence of automated collected data, more are more studies have used AFC and AVL data to monitor transit services. For example, \citet{wang2011bus} used AFC and AVL data to infer the OD demand patterns of the transit system in Transport for London (TfL). They calculated interchange time and the connecting bus route’s headway to evaluate the connectivity of bus networks. \citet{trepanier2009calculation} used smart card data to estimate the performance of transit systems by defining different metrics that can be directly extracted from the data. When combined with established evaluation processes, these measures enable operators to monitor their networks in greater detail. This includes assessing the performance of network supply and passenger service statistics at different spatial and temporal levels, including specific routes and bus stops. \citet{ma2014development} focused on creating a data-driven platform for online transit performance monitoring, utilizing data from AFC and AVL systems in Beijing. The platform employs data-mining techniques to estimate individual transit rider's origin and destination on flat-rate buses. This platform not only serves as a visualization tool for monitoring transit network performance but also contributes to data-driven transportation research and applications. \citet{nuzzolo2015dybus2} introduced a dynamic real-time transit simulation model. The model simulated travel behavior with real-time information, providing real-time predictions of on-board crowding. \citet{mo2020capacity} and \citet{mo2023ex} introduced a simulation-based optimization method to calibrate a transit network performance model. They used smart card data from the Mass Transit Railway network in Hong Kong as the case study, demonstrating the application of the proposed model for network performance monitoring. They also analyzed spatial-temporal crowding patterns in the system and evaluated different dispatching strategies. \citet{ma2019estimation} aimed to monitor rail transit system performance, particularly focusing on the number of times passengers are denied boarding. The research reviews existing methods and introduces a new approach for urban rail systems with both tap-in and tap-out fare collection systems. This method employs a mixture distribution framework with prior structural information, making it data-driven and not reliant on denied boarding observations or assumptions about access/egress time distributions.

\section{Methodology}\label{sec_method}

\subsection{Conceptual framework}\label{sec_concept}
As mentioned in Section \ref{intro}, in order to use a simulation model to monitor the transit network services, the two major tasks are:
\begin{itemize}
    \item Develop a transit simulation model and
    \item Calibrate the path choice parameters. 
\end{itemize}

\begin{figure}[htb] 
\centering
\includegraphics[width=0.8\linewidth]{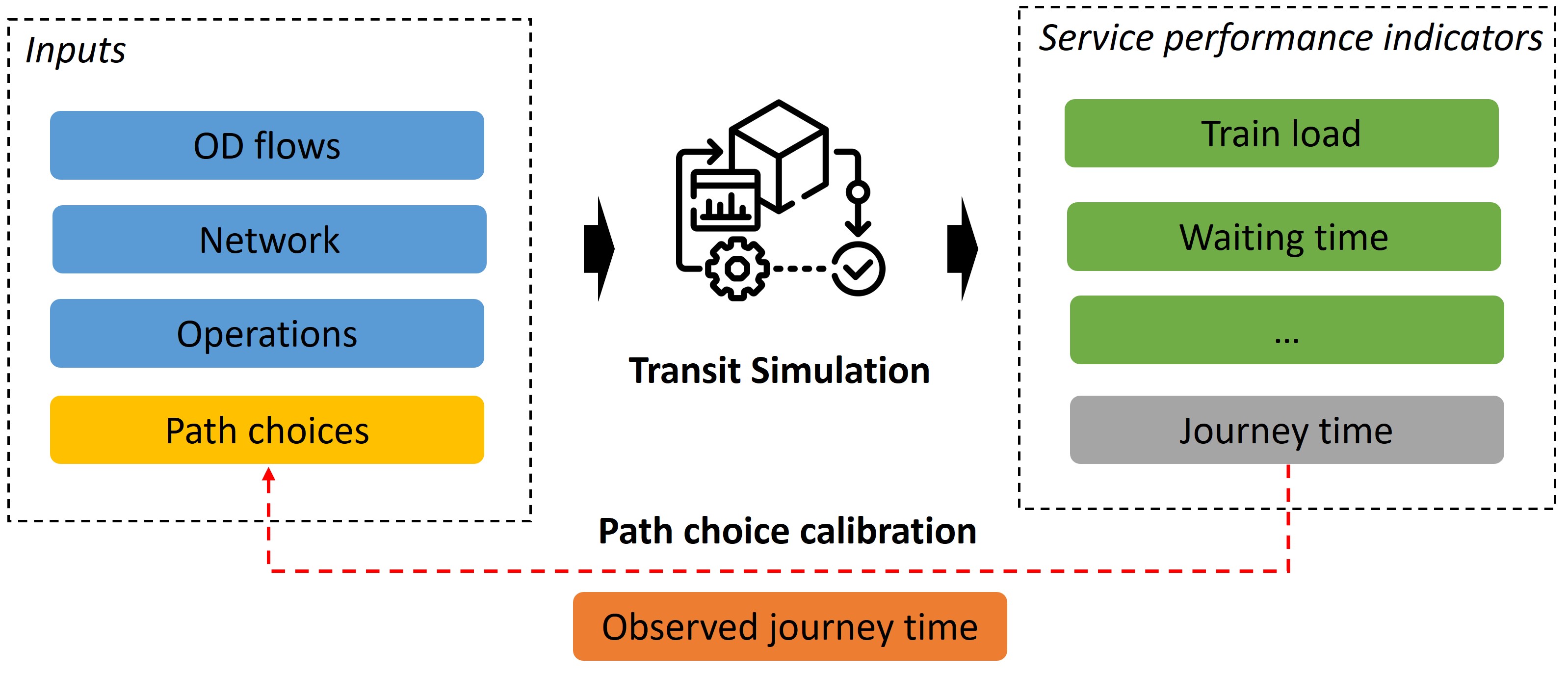}
\caption{Conceptual framework of the self-calibrated transit service monitoring}
\label{fig_concept}
\end{figure}

Consider a closed transit system (with both tap-in and tap-out records), passengers' journey time (i.e., tap-out time minus tap-in time) can be directly observed. Given a specific path, passengers' journey times are determined by the length of the path and the crowding of the path (i.e., they may not be able to board the first arrival train/bus if the vehicle is too crowded). Therefore, given different path choice patterns, passenger's journey time distribution can be different. As passenger journey time can be observed, this provides a way to calibrate path choices using the observed journey time: we wish to estimate a set of path choice parameters such that the model-output journey time distribution is close to the observed journey time distribution (as shown in Figure \ref{fig_concept}).

In the following sections, we first elaborate on how path choices are modeled (Section \ref{sec_path}) and the event-based transit simulation model (Section \ref{sec_sim}), then the calibration of path choices (Section \ref{sec_calibration}). 

\subsection{Path choice}\label{sec_path}
Path choice is usually modeled using the discrete choice framework, which assumes that decision makers maximize their utilities when making choices \citet{ben1985discrete}. The multinomial logit (MNL) model is a  typical example of the discrete choice model. For path choice problems, the C-logit model is often used. The C-logit is a variation of the MNL model which corrects for the fact that alternatives may not be independent due to path overlap. The C-logit incorporates an additional ``cost'' attribute, the commonality factor (denoted as $F_i$), in the utility \citep{cascetta1996modified}. Assume that the probability of any passenger choosing path $i$ of OD pair $w$ is:
\begin{align}
\mathbb{P}(i \in \mathcal{W}_w) = \frac{\exp{(\boldsymbol{\beta_X} \cdot \boldsymbol{X_{i}}+\beta_{F} \cdot F_i)}}{\sum_{j \in \mathcal{W}_w} \exp{(\boldsymbol{\beta_X} \cdot \boldsymbol{X_{j}} + \beta_{F} \cdot F_{j})}},
\label{eq_MNL}
\end{align}
where, $\boldsymbol{X_{i}}$ is the attribute vector of path $i$, such as in-vehicle time, number of transfers, etc. $\mathcal{W}_w$ is the set of all alternative paths for the OD pair $w$.  $\boldsymbol{\beta_X}$ and $\beta_{F}$ are the corresponding coefficients to be estimated. $F_i$ is the commonality factor of path $i$, defined as:
\begin{align}
F_i = \ln \sum_{j \in \mathcal{W}_w}(\frac{N_{i,j}}{N_i \cdot N_{j}})^\gamma
\label{eq_CF}
\end{align}
where $N_{i,j}$ is the number of common stations of paths $i$ and $j$. $N_i$ and $N_{j}$ are the number of stations for paths $i$ and $j$, respectively. $\gamma$ is a positive constant which is determined based on empirical studies \citep{li2014route}. In practice, the estimation of the route choice model parameters is based on survey data. In this study, we aim to estimate $\boldsymbol{\beta} = (\boldsymbol{\beta_X},\beta_{F})$ using the automated collected data (i.e., AFC and AVL).

\subsection{Transit simulation}\label{sec_sim}

\begin{figure}[htb] 
\centering
\includegraphics[width=0.8\linewidth]{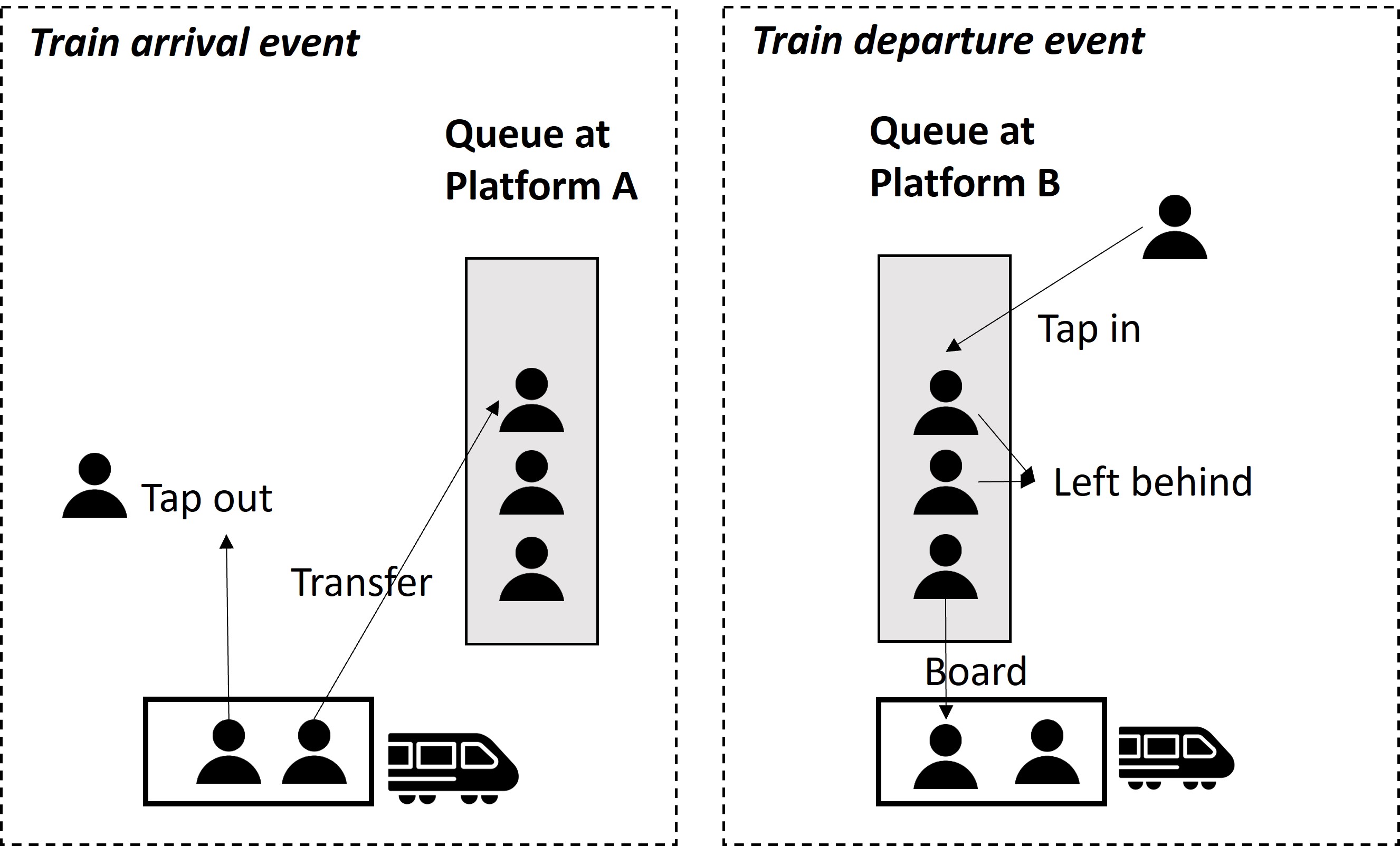}
\caption{Structure of the transit simulation model}
\label{fig_model}
\end{figure}

Figure \ref{fig_model} summarizes the main structure of the transit simulation model. Three objects are defined: trains, queues, and passengers. Trains are characterized by routes, runs, current locations, and capacities. Passengers are queued based on their arrival times. Three different types of passengers are represented: left-behind passengers who were denied boarding from previous trains, new tap-in passengers from outside the system, and new transfer passengers from other lines. The left-behind passengers are usually at the head of the queue. 

An event-based modeling framework is used to load the passengers onto the network. Two types of events are considered: train arrivals and train departures. The events are sorted by time and processed sequentially until all events are completed during the analysis period. When a train arrives at a station, the offloaded passengers either transfer or exit. Transfer passengers join the boarding queue. When a train departs a station, passengers are loaded on the train up to its available capacity based on a First Come First Serve (FCFS) principle.

\begin{algorithm} 
\caption{Transit Simulation} \label{alg_overall}
\begin{algorithmic}[1]
\State Get input demand and path choice probabilities
\State Generate event list based on train movement data (Section \ref{sec_event}). $K=$ total number of events
\State Set event counter $k=1$.
\Do
\State Assign passengers entering between event $k-1$ to $k$ to paths based on the input path choice probabilities (Section \ref{sec_pass})
\If{event $k$ is arrival}
\State Process train arrival event (Section \ref{sec_train_arr})
\Else
\State Process train departure event (Section \ref{sec_train_dep})
\EndIf
\State $k = k + 1$
\doWhile{$k \leq K$}
\State \Return Network performance indicators
\end{algorithmic}
\end{algorithm}

\subsubsection{Inputs of the simulation}\label{sec_input}
The inputs of the simulation model include three components. The first one is the path choice parameter $\boldsymbol{\beta}$. The second is the tap-in time of each passenger in the system, referred to as OD entry demand. Denote all OD entry demand as $\boldsymbol{D}$. The last part is operation parameters, such as timetable, train configurations, and train network layout. Denote all operation parameters as $\boldsymbol{C}$. The simulation process can be described as
\begin{align}
    \text{Transit simulation}(\boldsymbol{\beta};\;\boldsymbol{D}, \boldsymbol{C}) 
\end{align}

\subsubsection{Event generation}\label{sec_event}
\begin{figure}[htb] 
\centering
\includegraphics[width=0.8\linewidth]{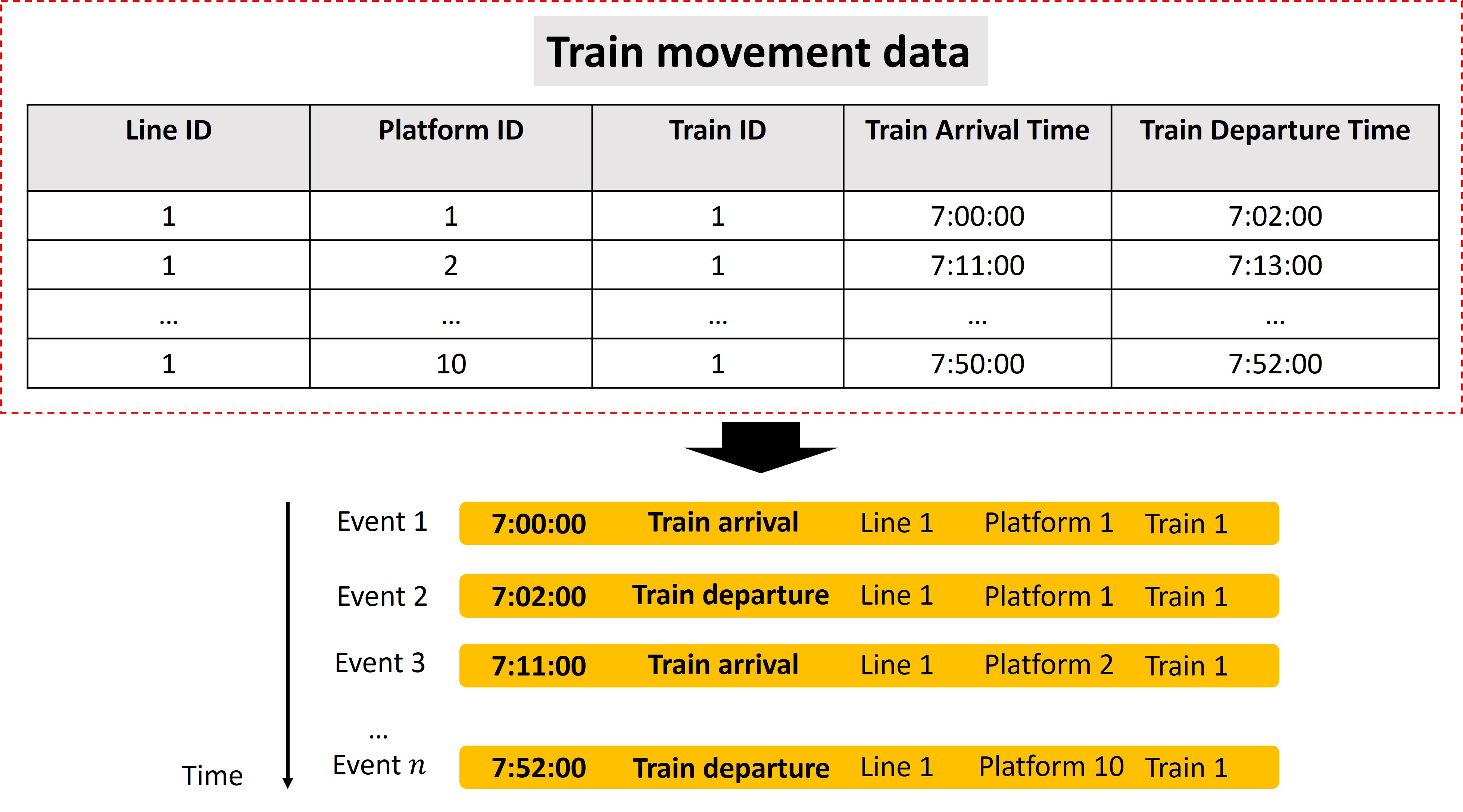}
\caption{Illustration of event generation from train movement data}
\label{fig_event}
\end{figure}

The event lists (arrivals and departures) of the simulation model are generated according to the actual train movement data from AVL (or timetable). Each event contains a train ID, occurrence time, and location (platform). Figure \ref{fig_event} presents an example of how a series of events are generated from a sample timetable. 

\subsubsection{Passenger path assignment}\label{sec_pass}
At the beginning of event $k$, we first obtain all passengers entering between event $k-1$ and $k$. Then each passenger is randomly assigned to a path based on the corresponding path choice probability estimated from a path choice model (Section \ref{sec_path}). Given the assigned path, the passengers are put into the queue at the first station of the path according to their arrival time at the platform.

\subsubsection{Train arrivals}\label{sec_train_arr}
For an arrival event, the train offloads passengers who reach their destination or need to transfer to the station. After offloading, train load and in-vehicle passengers of the train are updated. For passengers who reach their destinations, we remove them from the system and record their tap-out times. Their tap-out times are calculated based on the offloading time plus their egress time. Passengers who transfer at the station will be added to the waiting queue at the next boarding station. Their sequences in the queue are based on their arrival times at the next platform, which are calculated based on the offloading time plus transfer time. Note that the egress time of a station and transfer time between stations are both pre-determined parameters according to the layouts of stations.

\subsubsection{Train departures}\label{sec_train_dep}
For departure events, we first add new tap-in passengers (i.e., passengers who arrive at the platform between the last train departure and now) into the queue at the platform. Their sequences at the queue are based on their arrival times, which are calculated as tap-in time plus access time. Access time is also a pre-determined parameter. Then, passengers board the train according to an FCFS discipline until the train reaches its capacity. Passengers who cannot board are left behind and wait in the queue for the next train. The states of the train and the waiting queue are then updated accordingly.

\subsection{Calibration of path choices}\label{sec_calibration}
As discussed in Section \ref{sec_concept}, passengers' path choices will impact their journey time in the network, where the ground truth values can be directly observed in the AFC data. Therefore, the path choice calibration can formulate an optimization problem. The objective function is to minimize the difference between observed and model-output passenger journey times, with path choice parameters $\boldsymbol{\beta}$ as decision variables. 

For each passenger $n$ in the AFC data records, denote its observed journey time as $\hat{J}_n$ and the model-output journey time as $J_n$. Directly minimizing the difference between $J_n$ and $\hat{J}_n$ is too granular considering other random factors in the systems (like walking speed, and other activities in the transit system like shopping). Therefore, we need to consider a more aggregated way to model the journey time differences of all passengers. Following the same idea in \citet{mo2021calibrating}, we propose two aggregated indicators: 1) OD exit flows and 2) Kullback-Leibler (KL) divergence. 

Let us consider discrete time intervals $t\in\mathcal{T}$ in the system, where each time interval has the same length $\tau$ (e.g., $\tau = $15 min), $\mathcal{T}$ is the set of all time intervals. The OD exit flow of origin $o$ and destination $d$ at time interval $t$ is defined as the total number of passengers exiting at station $d$ at time interval $t$ with origin $o$ and destination $d$. Denote the observed OD exit flow from AFC data as $\hat{q}_{o,d,t}$ and the model-output values as ${q}_{o,d,t}$. It is worth noting that ${q}_{o,d,t}$ captures some of the journey time information because passengers with the same OD pair may contribute to different $q_{o,d,t}$ because of the differences in journey time. 

Another indicator is the KL divergence, which describes the distribution differences. Let $J_{o,d,t}$ be the journey time of passengers departing at time $t$ and OD pair $(o,d)$. Denote $\mathbb{P}_{o,d,t}[J_{o,d,t}]$ as the model-output journey time distribution and  $\hat{\mathbb{P}}_{o,d,t}[J_{o,d,t}]$ as the observed one. Then the KL divergence is
\begin{align}
D^{\text{KL}}_{o,d,t}\left(\mathbb{P}_{o,d,t}[J_{o,d,t}]\;||\;\hat{\mathbb{P}}_{o,d,t}[J_{o,d,t}]\right) = \int_{x} \mathbb{P}_{o,d,t}[x] \cdot \log{\left(\frac{\mathbb{P}_{o,d,t}[x]}{\hat{\mathbb{P}}_{o,d,t}[x]}\right)}\text{d}x
\label{eq_DL}
\end{align}
Note that in the case study, we consider a discrete distribution of $J_{o,d,t}$, which simplifies the integration to summation for KL divergence calculation.

Combining the KL divergence and OD exit flow, the final optimization model is formulated as:
\begin{subequations}\label{eq_opt}
\begin{align}
    \min_{\boldsymbol{\beta}} &\quad Z(\boldsymbol{\beta}) = \sum_{(o,d)\in\mathcal{W},t\in\mathcal{T}} \left(q_{o,d,t} - \hat{q}_{o,d,t} \right)^2 + \eta \cdot \sum_{(o,d,t)\in\mathcal{WT}^*} D^{\text{KL}}_{o,d,t}\left(\mathbb{P}_{o,d,t}[J_{o,d,t}]\;||\;\hat{\mathbb{P}}_{o,d,t}[J_{o,d,t}]\right) \\
    \text{s.t. } & \mathbb{P}_{o,d,t}[J_{o,d,t}]= \text{Transit simulation}(\boldsymbol{\beta};\;\boldsymbol{D}, \boldsymbol{C}) \quad \forall (o,d,t)\in\mathcal{WT}^*\\
    & q_{o,d,t} = \text{Transit simulation}(\boldsymbol{\beta};\;\boldsymbol{D}, \boldsymbol{C}) \quad \forall (o,d)\in\mathcal{W},t\in\mathcal{T}
\end{align}
\end{subequations}
where $\eta$ is the trade-off parameter for two objective function terms. $\mathcal{W}$ is the set of all OD pairs. $(o,d,t)\in\mathcal{WT}^*$ is a subset of of $\mathcal{W}\times\mathcal{T}$ with $\hat{q}_{o,d,t}$ greater than or equal to a pre-determined threshold $Q^\text{KL}$. Mathematically:
\begin{align}
    \mathcal{WT}^* = \{(o,d,t):\; (o,d)\in\mathcal{W}, t\in\mathcal{T}, \hat{q}_{o,d,t} \geq Q^\text{KL}\}
\end{align}
The threshold $Q^\text{KL}$ is used to ensure the calculation of journey time distribution has enough samples. 

Eq. \ref{eq_opt} has a non-analytical constraint with a simulation model. It can be solved through simulation-based optimization. In this study, we adopt the constrained optimization with response surfaces (CORS) method in \citet{knysh2016blackbox} to solve the problem. The main idea is to approximate the original objective function $Z(\boldsymbol{\beta})$ with a surrogate model $\hat{Z}(\boldsymbol{\beta})$ (e.g., radial basis functions) such that $\hat{Z}(\boldsymbol{\beta})$ has an analytical form. Then we perform typical convex optimization methods on $\hat{Z}(\boldsymbol{\beta})$ to solve for the $\boldsymbol{\beta}$. With new $\boldsymbol{\beta}$, we can update the fitted model $\hat{Z}(\boldsymbol{\beta})$ and solve it again. Repeat this process until we reach maximum iterations or other convergence criteria.

\section{Case study}\label{ref_case}
\subsection{Hong Kong Mass Transit Railway}
\begin{figure}[htb] 
\centering
\includegraphics[width=0.8\linewidth]{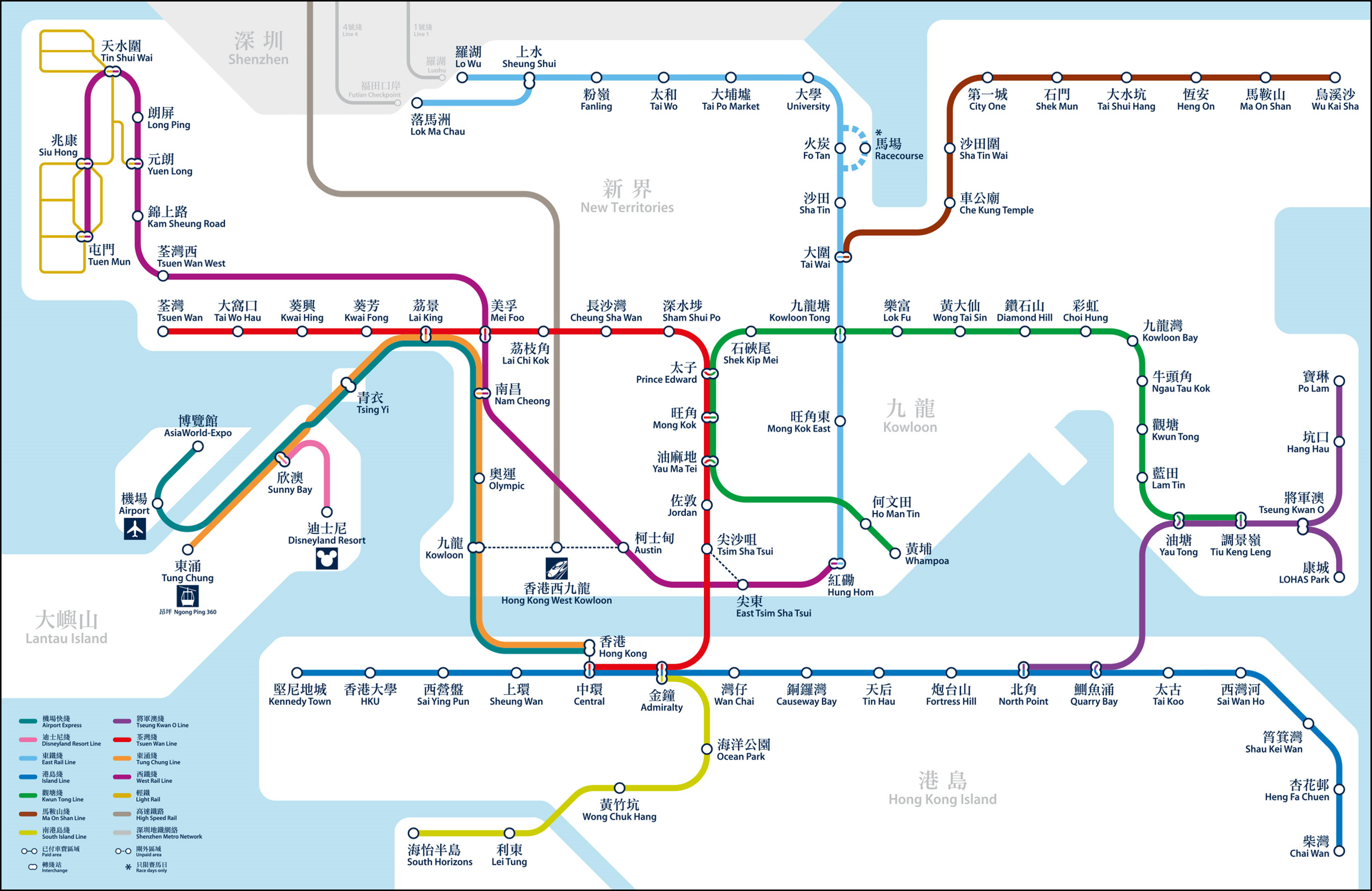}
\caption{Hong Kong MTR system}
\label{fig_mtr}
\end{figure}

The proposed transit service monitoring model is verified through a case study of the Hong Kong MTR system (Figure \ref{fig_mtr}). The MTR system covers urbanized regions in Hong Kong, including Hong Kong Island, Kowloon, and the New Territories. It comprises 11 lines spanning a total of 218.2 km (135.6 miles) of rail, with 159 stations, including 91 heavy rail stations and 68 light rail stops. On an average weekday, the network facilitates over 5 million trips.

\subsection{Synthetic data generation}
In the case study, we will generate the synthetic tap-in and tap-out data for model validation. Specifically, we assume all passengers follow an ``actual'' path choice behavior quantified by a discrete choice model with a synthetic vector of $\boldsymbol{\beta}^*$ (Table \ref{tab_beta}). The values of $\boldsymbol{\beta}^*$ are obtained from \citet{mo2020network}. The main explanatory variables are the total in-vehicle time, relative walk time, and the number of transfers. The relative walk time is defined as the total walk time (access + transfers + egress) divided by the map distance of the path. Parameters for these attributes are all negative, implying that routes with high in-vehicle, walk, and transfer times are less likely to be chosen by passengers. Then, we input the actual passenger tap-in time and $\boldsymbol{\beta}^*$ to the transit simulation model to generate the passenger tap-out time. Collecting all passengers' tap-in and tap-out times, we obtain the synthetic AFC data. Train operation is assumed to follow the MTR schedules. 

\begin{figure}[htb] 
\centering
\includegraphics[width=0.7\linewidth]{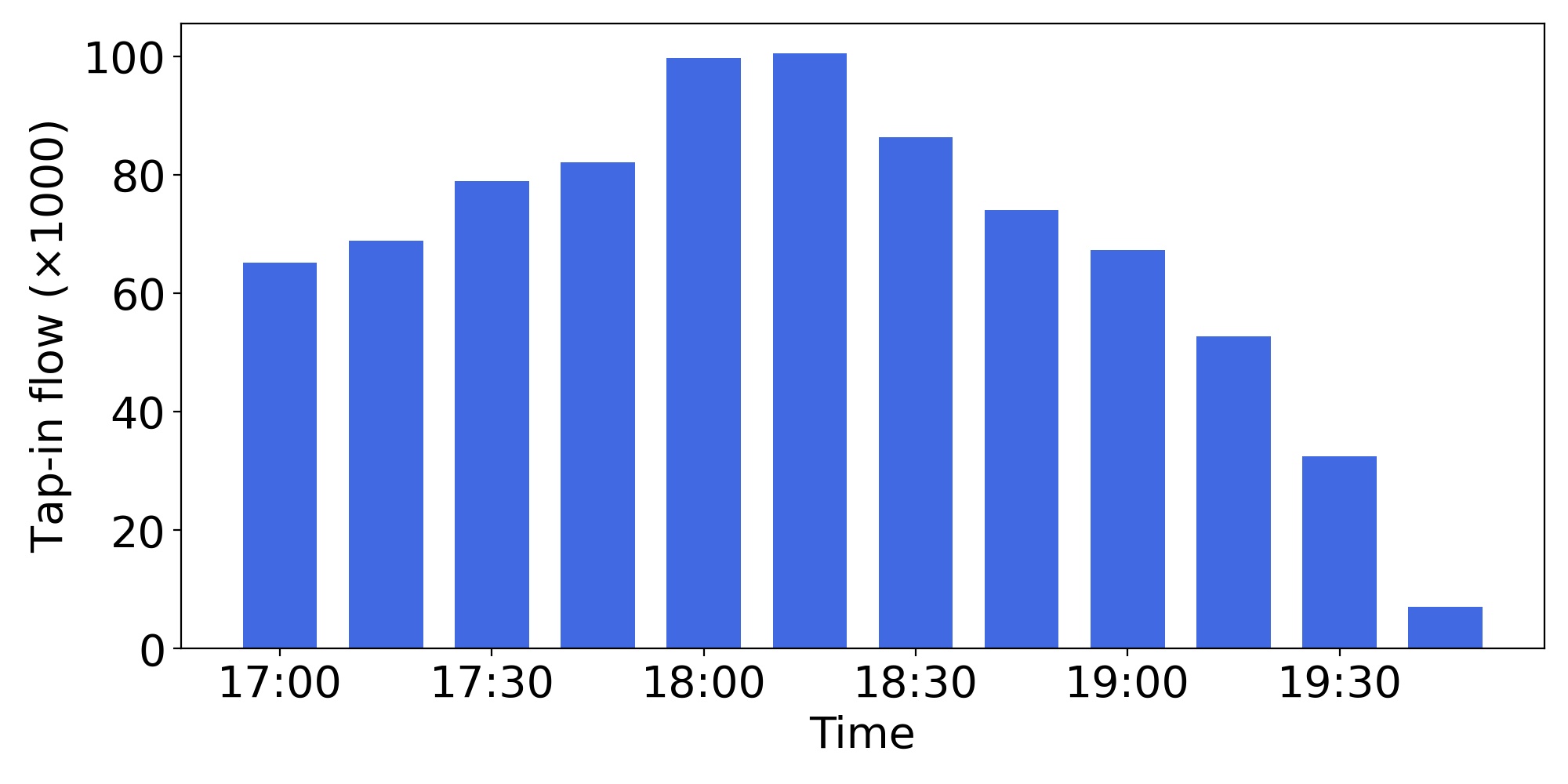}
\caption{Tap-in flow distribution}
\label{fig_tapin}
\end{figure}

Passenger tap-in flows used in the study are shown in Figure \ref{fig_tapin}. The total study period is 17:00 - 20:00, the evening peak hours of the system. The first 1 hour is set as the warm-up time and the last 1 hour is set as the cool-down time. Only the hour for 18:00 - 19:00 is used for the path choice estimation. During the three hours, more than 800 thousand passengers are using the system. The length of each time interval is set as 15 minutes. 

\subsection{Parameter settings}
The parameter for the commonality factor is set as $\gamma = 5$. The trade-off parameter $\eta$ is set as 600 to account for the scalar and importance difference between two objective terms. Passenger waiting time is assumed to be log-normally distributed. We assume the average speed is 1.5 m/s. The distance between stations and from gates to platforms is obtained from \citep{mtr_layout}. The minimal number of passengers to calculate journey time distribution is set as $Q^\text{KL} = 50$. The maximum iteration for the optimization is set as 100.

\subsection{Results}
\begin{figure}[htb] 
\centering
\includegraphics[width=0.8\linewidth]{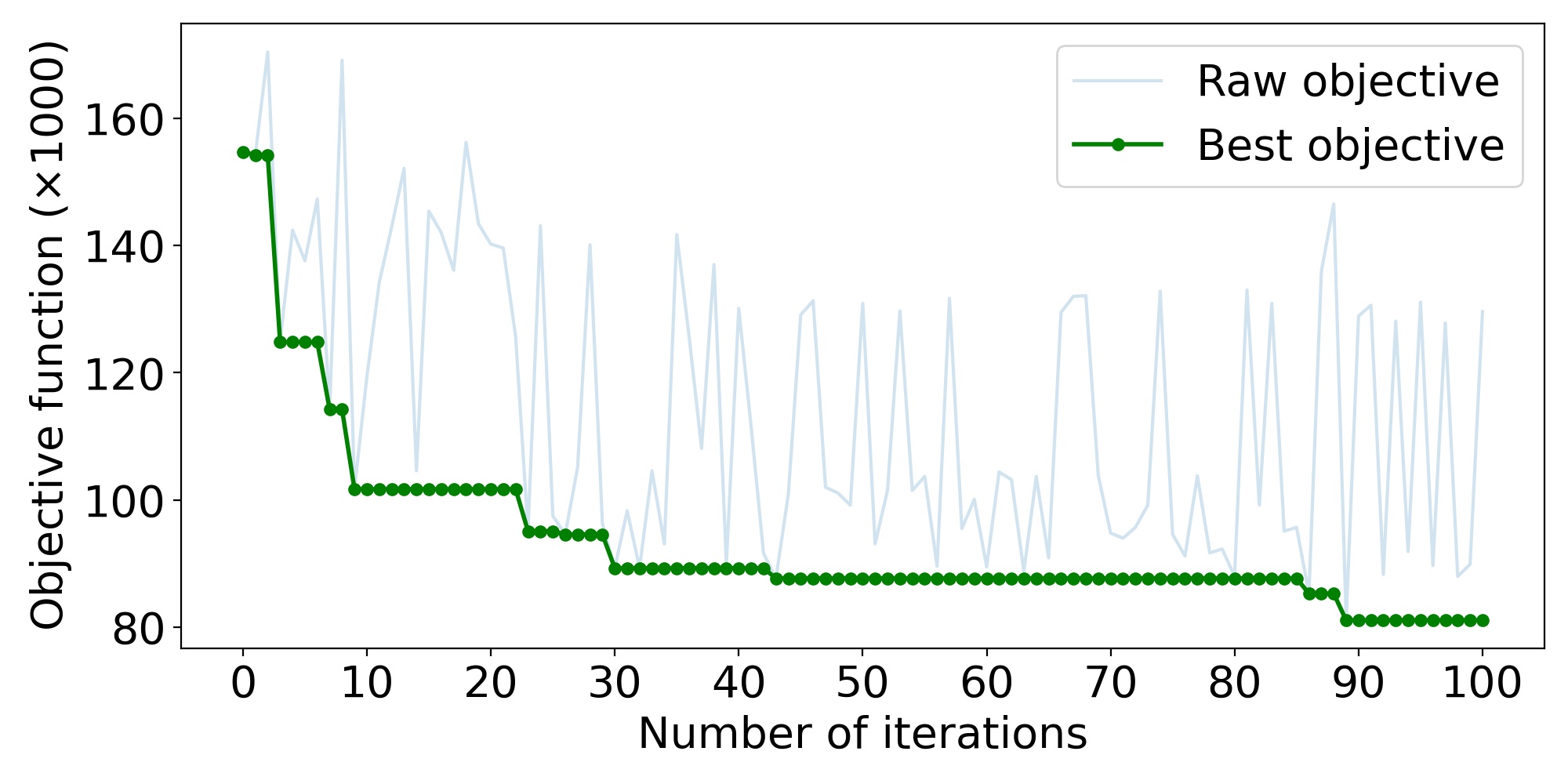}
\caption{Objective function}
\label{fig_obj}
\end{figure}

The objective function during the optimization process is shown in Figure \ref{fig_obj}. Since the CORS method is a heuristic, the raw objective function at each iteration is not guaranteed to be decreased with iterations. After 100 iterations, we get a relatively good result. The estimated $\boldsymbol{\beta}$ are shown in Table \ref{tab_beta}. The results are close to the ``true'' parameters, which validate the effectiveness of the proposed model.

\begin{table}[htb]
\caption{Optimal $\beta$ Estimation Results}
\centering
\begin{tabular}{@{}ccccc@{}}
\toprule
          & In-vehicle time & Relative walk time & Number of transfers & Commonality factor \\ \midrule
``True'' ($\boldsymbol{\beta}^*$)  & -0.147         &  -1.271       &     -0.573         &    -3.679        \\
Estimated &  -0.183            &    -1.666                &      -0.715               &  -4.852                  \\ \bottomrule
\end{tabular}
\label{tab_beta}
\end{table}

\begin{figure}[htb]
\centering
\subfloat[18:00 - 18:30]{\includegraphics[width=1.0\textwidth]{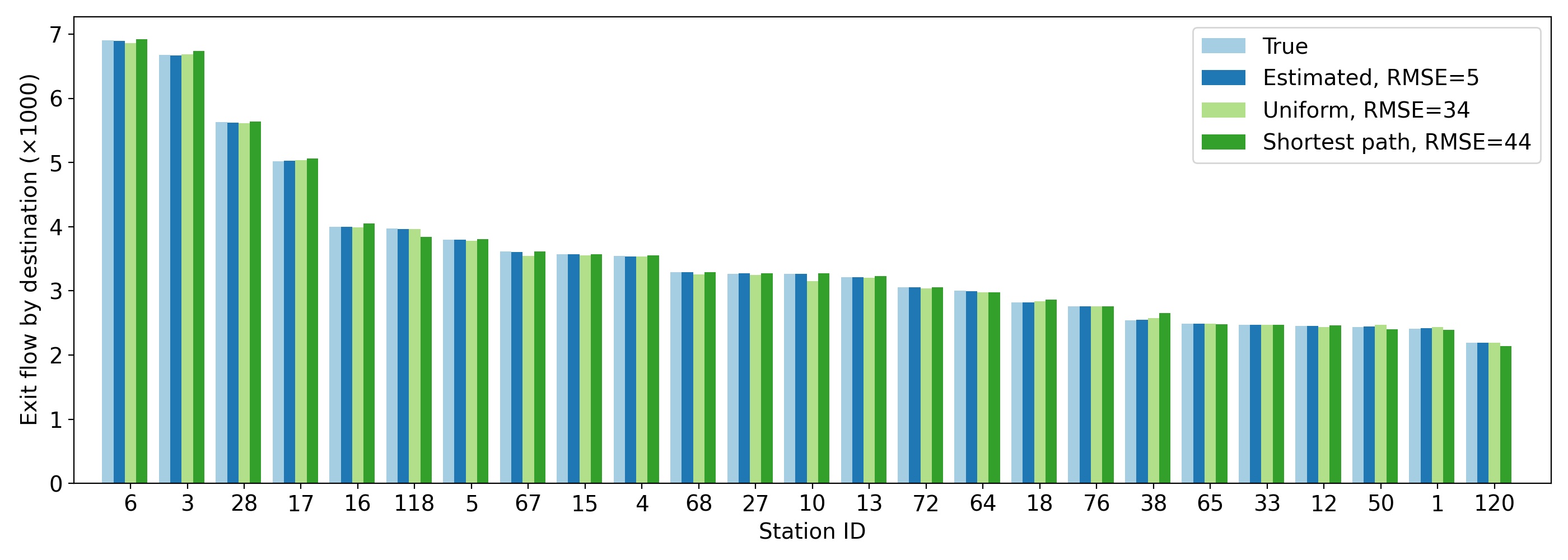}\label{fig_time1}}
\hfil
\subfloat[18:30 - 19:00]{\includegraphics[width=1.0\textwidth]{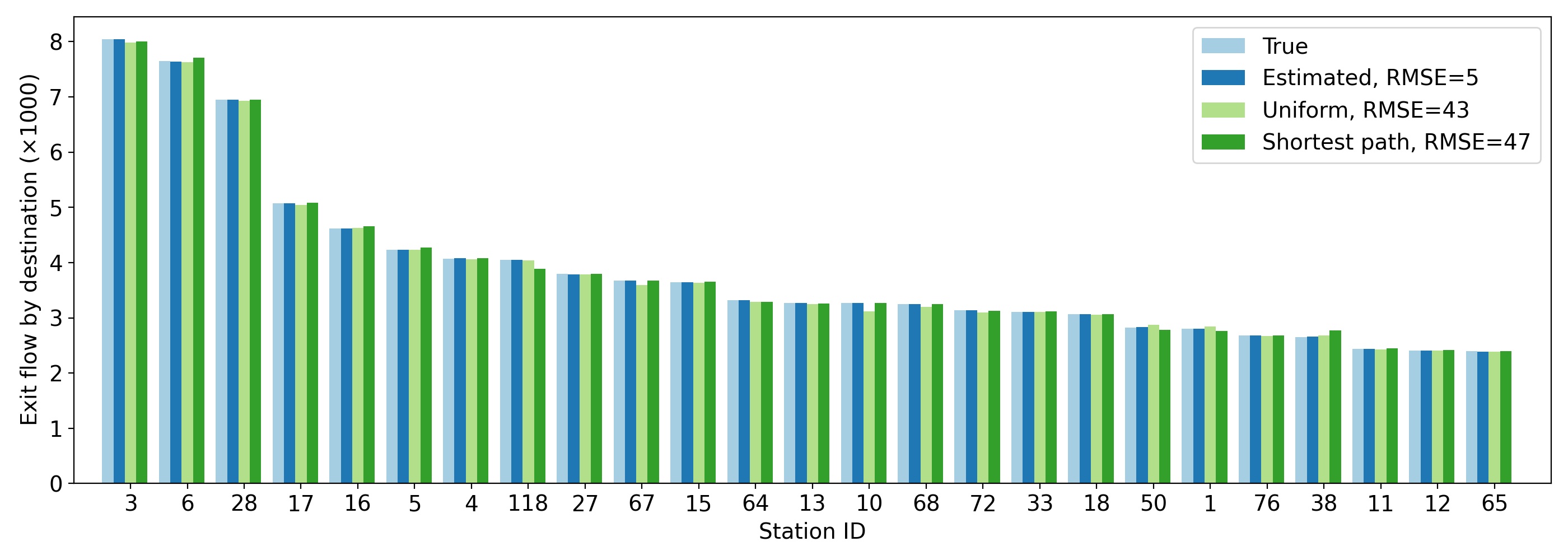}\label{fig_time2}}
\caption{Comparison of OD exit flows}
\label{fig_compare_od}
\end{figure}

We also compare the estimated OD exit flow with the ``actual'' ones. The exit flows are aggregated by time and origins. For comparison purposes, we compare the proposed model with 2 benchmark models:
\begin{itemize}
    \item \textbf{Uniform path choices}: Uniform path choices assume all passengers have equal probabilities of choosing all available paths. This is equivalent to assuming $\boldsymbol{\beta} = 0$.
    \item \textbf{Shortest path choices}: In the shortest path model, we assume all passengers will choose the shortest path among the OD pair, which is equivalent to $\beta_{\text{In-veh time}} = -\infty$ all other parameters are 0.
\end{itemize}
The OD exit flows are obtained by inputting different path choice parameters into the simulation model. Figure \ref{fig_compare_od} shows the comparison in two different time periods (18:00 - 18:30 and 18:30 - 19:00). The OD exit flows with optimization-based $\boldsymbol{\beta}$ are closer to the ``true'' OD exit flows (with ``true''  $\boldsymbol{\beta}$ as input) for both time periods with much lower root mean square errors (RMSE). This shows that the proposed model can output more accurate transit performance indicators compared to benchmark methods. 

\section{Conclusion and Discussion}\label{sec_con_dis}
This paper proposes a self-calibrated transit service monitoring framework that aims to obtain the performance of a transit system using automated collected data. We first introduce an event-based transit simulation model to simulate passenger travel behavior in a transit system. To estimate passenger path choices, we assume the path choices can be modeled using a C-logit model, and propose a simulation-based optimization model to estimate the path choice parameters based on AFC and AVL data. The Hong Kong MTR system is used as the case study. Results show that the model can well estimate the path choice behavior of passengers in the system. The output passenger exit flows are closer to the actual one compared to the two benchmark models (shortest path and uniform path choice). 

The proposed model has the potential to be implemented in the real-world. It eliminates the need for additional equipment (such as cameras) to collect network performance indicators. Instead, all data inputs can be obtained from automated data collection systems. If all AFC and AVL data can be obtained in real-time, the model can be run in real-time to simulate the network performance. Besides, the model can also be used to evaluate the historical performance with automatic daily calibration.

\bibliography{mybibfile}

\end{document}